\documentclass[letterpaper, 12pt]{article}
\usepackage{amssymb,amsmath}
\usepackage{epsf}
\usepackage{times,bm,mathrsfs}  % Fonts
\usepackage{amsmath}
\usepackage{amssymb}
\usepackage{amsfonts}
\usepackage{mathrsfs}
\usepackage{bbm}
\usepackage{color}
\usepackage{graphicx}
\usepackage{amscd}
\usepackage{epsfig}

\usepackage{setspace}
\definecolor{Red}{rgb}{1,0,0}
\definecolor{Blue}{rgb}{0,0,1}
\definecolor{Olive}{rgb}{0.41,0.55,0.13}
\definecolor{Green}{rgb}{0,1,0}
\definecolor{MGreen}{rgb}{0,0.8,0}
\definecolor{DGreen}{rgb}{0,0.55,0}
\definecolor{Yellow}{rgb}{1,1,0}
\definecolor{Cyan}{rgb}{0,1,1}
\definecolor{Magenta}{rgb}{1,0,1}
\definecolor{Orange}{rgb}{1,.5,0}
\definecolor{Violet}{rgb}{.5,0,.5}
\definecolor{Purple}{rgb}{.75,0,.25}
\definecolor{Brown}{rgb}{.75,.5,.25}
\definecolor{Grey}{rgb}{.5,.5,.5}
\definecolor{Black}{rgb}{0,0,0}

\def\path{{\tt path}}

\newcommand{\fcal}{\mathcal{F}}

\newcommand{\real}{\mathbb{R}}

\newcommand{\eps}{\varepsilon}

\newcommand{\bdm}{\begin{displaymath}}
\newcommand{\edm}{\end{displaymath}}
\newcommand{\bea}{\begin{eqnarray*}}
\newcommand{\eea}{\end{eqnarray*}}
\newcommand{\bean}{\begin{eqnarray}}
\newcommand{\eean}{\end{eqnarray}}

\newcommand{\expec}{\mathbb{E}}

\newcommand{\poly}{\mathrm{poly}}

\newcommand{\sm}{\setminus}

\newcommand{\bfaa}{\mathbf{A}}
\newcommand{\bfbb}{\mathbf{B}}
\newcommand{\bfcc}{\mathbf{C}}
\newcommand{\bfdd}{\mathbf{D}}
\newcommand{\bfss}{\mathbf{S}}
\newcommand{\bftt}{\mathbf{T}}

\newcommand{\qcal}{\mathcal{Q}}

\newcommand{\NP}{\mathrm{NP}}

\newtheorem{theorem}{Theorem}

\newtheorem{corollary}{Corollary}
\newtheorem{definition}{Definition}

\newtheorem{lemma}{Lemma}
\newtheorem{conjecture}{Conjecture}

\newtheorem{assumption}{Assumption}

\newenvironment{proof}{\noindent{\textbf{Proof:}}}{$\blacksquare$\vskip\belowdisplayskip}

\begin{document}

\title{
%\vspace{-3cm}Local to Global Submodularity in Social Networks
\vspace{-3cm}Submodularity of Infuence in Social Networks: From Local to Global
\thanks{
 E.M. is supported by an Alfred Sloan fellowship in
Mathematics and by NSF grants DMS-0528488, and DMS-0548249 (CAREER) and by 
ONR grant N0014-07-1-05-06.
}
}
\author{
Elchanan Mossel\footnote{Statistics and Computer Science, UC Berkeley
and Mathematics and Computer Science, Weizmann Institute.} and
Sebastien Roch\footnote{Mathematics, UCLA. Work done at UC Berkeley.}
}
\maketitle

\begin{abstract}
{
%We prove and extend a conjecture of Kempe,
%Kleinberg, and Tardos (KKT) on the spread of influence in social networks. 
%
Social networks are often represented as directed graphs 
where the nodes are individuals and the edges indicate a form of social relationship.
A simple way to model the diffusion of ideas, innovative behavior, or ``word-of-mouth'' effects
on such a graph is to consider an increasing process of ``infected'' (or active) nodes:
each node becomes infected once an activation function of the set of its
infected neighbors crosses a certain threshold value. 
Such a model was introduced
by Kempe,
Kleinberg, and Tardos (KKT) in~\cite{KeKlTa:03,KeKlTa:05} where the authors
also impose several natural assumptions: the threshold
values are random and the activation functions
are monotone and submodular. The monotonicity condition
indicates that a node is more likely to become active if more of its neighbors
are active, while the submodularity condition 
indicates that the marginal effect
of each neighbor is decreasing when the set of active
neighbors increases.

For an initial set of active nodes $S$, let 
$\sigma(S)$ 
denote the expected number of active nodes at termination. 
Here we prove a conjecture of KKT: we show that the function $\sigma(S)$ is
submodular under the assumptions above. 
We prove the same result for the expected value 
of any monotone, submodular function of the set of active nodes at termination.   
Roughly,
our results demonstrate that ``local'' submodularity is preserved
``globally'' under this diffusion process. This is of natural computational interest, as
many optimization problems have good approximation algorithms for 
submodular functions. 
%In particular, our results coupled with an argument in~\cite{KeKlTa:03} imply that a greedy
%algorithm gives an $(1-1/e-\eps)$-approximation algorithm for maximizing 
%$\sigma(S)$ among all sets $S$ of a given size. This
%\vfill\eject
%\noindent result has important
%practical implications for many social network analysis problems, notably
%viral marketing.
}
\end{abstract}

%\category{G.3}{Probability and Statistics}{Stochastic Processes}
%A category including the fourth, optional field follows...
%\category{D.2.8}{Software Engineering}{Metrics}[complexity measures, performance measures]

%\terms{Algorithms, Economics, Theory}

% EL: Changed here

%\keywords{Viral marketing, social networks, submodularity, coupling}

\section{Introduction}
{\bf Social Networks.}
{
In recent years, diffusion processes on social networks have been the focus
of intense study. While traditionally such processes have been
of major interest in epidemiology where they model the spread of diseases and immunization~\cite{Morris:04,Liggett:85,Liggett:01,Durrett:88,BeBoChSa:05,DuJu:u},
much of the recent interest has resulted 
from applications in sociology, 
economics, and engineering~\cite{BrownReinegen:87,
AsRoLeVe:01, GoLiMu:01a,GoLiMu:01b,
DomingosRichardson:01,DomingosRichardson:02,
KeKlTa:03,KeKlTa:05}.

In computer science, a strong motivation for analyzing diffusion processes
has recently emanated 
from the study of viral marketing strategies in data mining, 
where various novel algorithmic problems
have been considered~\cite{DomingosRichardson:01,DomingosRichardson:02,
KeKlTa:03,KeKlTa:05}. Roughly speaking, {\em viral marketing}---unlike conventional mar\-keting---takes 
into account the ``network value'' of potential customers, that is, it seeks to target
a set of individuals whose influence
on their social network through word-of-mouth effects is high.
%For more background on viral marketing, see~\cite{DomingosRichardson:01,DomingosRichardson:02,
%KeKlTa:03,KeKlTa:05}.

Commonly-used heuristics to identify influential nodes in social networks
include picking individuals of high degree---so-called degree centrality heuristics---or
picking individuals with short average distance to the rest of the network---so-called
distance centrality heuristics~\cite{WassermanFaust:94}.
Typically these heuristics provide no guarantees.
Here we prove a structural 
conjecture of Kempe, Kleinberg, and Tardos (KKT)~\cite{KeKlTa:03, KeKlTa:05}
which in turn implies that a natural greedy algorithm has a good
performance guarantee.
The conjecture can be roughly stated as follows: if a diffusion model
is locally submodular, that is, the influence on an individual by its neighbors
in the network has ``diminishing returns,'' then the process is globally
submodular (on average). This is relevant in this context because, under the submodularity property,
optimization problems---such as the viral marketing problem---are known to have good approximation 
algorithms~\cite{NemhauserWolsey:88}. 
% EL: Changed that have to to have
In particular, in~\cite{KeKlTa:03}, greedy algorithms based on the above conjecture 
were shown to achieve significantly better performances in
practice than widely-used network analysis heuristics.

{\bf General Threshold Model.} In~\cite{KeKlTa:03}, KKT introduced
the \emph{general threshold model}, a broad
generalization of a variety of natural diffusion models on networks, including
the influential \emph{linear threshold model} of Granovetter in sociology~\cite{Granovetter:78}. 
Given an initial set of infected or active individuals
on a network, the process grows in the following way. 
(See Section~\ref{sec:model} for a formal description.)
Each individual, say $v$, has an activation function, which  
measures the effect of its neighbors on $v$, as well as a 
threshold value. At any time, if the set of previously infected neighbors
of $v$ is such that its activation function crosses its threshold value,
then $v$ becomes infected.
This process is {\em progressive}---an
active node stays active forever. KKT consider the following natural assumptions:
\begin{itemize}
\item[-]
The {\em threshold values} are {\em random}. This is to account 
for our lack of knowledge of the exact threshold values.
KKT actually assume that the thresholds are \emph{uniformly} random. 
Note however that, given that any distribution can be generated from
a uniform random variable, it follows that, by 
appropriately modifying the activation functions, the threshold
values can \emph{effectively} have any distribution.
(See Section~\ref{sec:conclusion}.)
\item[-]
The {\em activation functions} are {\em monotone increasing}. That is, 
a node is more likely to become infected if a larger set of its neighbors is infected.
\item[-]
The {\em activation functions} are {\em submodular}. This corresponds to the fact that
the marginal effect of each neighbor of $v$ decreases as the set of
active nodes increases.
\end{itemize}

{\bf The Influence Maximization Problem.}
Since the diffusion process defined above is increasing, it terminates after a finite number of steps. 
For a given initial set
of active nodes $S$ we define $\sigma(S)$ to be the expected size of
the set of active nodes at the end of the process. 
In the {\em Influence Maximization Problem}, we aim to find a set $S$ of a
fixed size maximizing $\sigma(S)$. 

The Influence Maximization Problem is a natural problem to consider in
the context of viral marketing. Given a social network, 
it is desired to find a small set of ``target'' individuals 
so as to maximize the number of customers who will eventually purchase a
product following the effects of ``word-of-mouth''~\cite{DomingosRichardson:01,DomingosRichardson:02}. 
The same problem may also be of interest in epidemiology where finding
the set $S$ of a fixed size maximizing $\sigma(S)$ is a natural problem
both in terms of bounding the spread of a disease and in terms of
maximizing the effect of immunization. 

In~\cite{KeKlTa:03} it was shown that the Influence Maximization Problem
is $\NP$-hard to approximate within a
factor $1-1/e+\eps$ for all $\eps > 0$ and that the problem is in fact
$n^{1-\eps}$ hard to approximate without the submodularity
condition. (See, e.g.,~\cite{GareyJohnson:79,ACGKMP:99} for background
on $\NP$-completeness and hardness of approximation.)
On the other hand, it was shown in~\cite{KeKlTa:05} that for all $\eps > 0$ it is
possible to find a set $S$ of fixed size that is a $(1-1/e-\eps)$-approximation 
of the maximum in random polynomial time if the set function $\sigma$ is {\em itself}
submodular. This leads to the following conjecture. 

\begin{conjecture}[\cite{KeKlTa:03,KeKlTa:05}] \label{conj}
The function $\sigma$ is submodular.
\end{conjecture} 

While the result of~\cite{KeKlTa:03,KeKlTa:05} showed that $\sigma$ is
submodular in special cases and related models (see below), the general case was open
prior to our work.  
In this paper we prove Conjecture~\ref{conj} and extend it to the
case where $\sigma(S)$ is the expected value of any monotone, submodular function
of the final active set. This gives a $(1-1/e-\eps)$-approximation algorithm 
for finding a set $S$ of fixed size maximizing $\sigma(S)$.

\subsection{The Model}\label{sec:model}
%EL: Added The

In this section, we define formally the {\em general threshold model}.
\begin{definition}[Social Network]
A {\em social network} is given by:
\begin{itemize}
\item[-]
A {\em ground set} $V$ with $|V| = n$ 
\item[-]
A collection of {\em activation functions} 
$\fcal = (f_v)_{v\in V}$, where $f_v : 2^V \to [0,1]$ is a $[0,1]$-valued set function
on $V$. 
\end{itemize}
\end{definition}
Typically, we think of $V$ as the individuals of a social network 
$G = (V,E)$ where each $f_v$ measures the effect of $v$'s neighbors
$N(v)$ on $v$. In particular $f_v$ depends only on neighbors $N(v)$ 
affecting $v$, so $f_v(S) = f_v(N(v)\cap S)$ for all $S$. 
However, the specification of the graph will not be needed below. 

\begin{definition}[Monotonicity]
The function\\ $f : 2^V \to \real$ is {\em monotone} if 
\[
f_v(S) \leq f_v(T)
\]
for all 
\[
S\subseteq T \subseteq V.
\]
We say that a collection $\fcal$ of functions from $2^V$ to $\real$
is monotone if all its elements are monotone functions.
\end{definition}
\begin{definition}[Submodularity]
The function\\ $f: 2^V \to \real$ is {\em submodular}
if for all $S,T \subseteq V$
\begin{equation*}
f(S) + f(T) \geq f(S\cap T) + f(S\cup T).
\end{equation*}
We say that a collection $\fcal$ of functions from $2^V$ to $\real$
is submodular if all its elements are submodular functions.
\end{definition}
The monotonicity condition corresponds to the fact that the effect of
a larger set on $v$ is stronger than the effect of a smaller set. 
The submodularity condition is equivalent to the fact that if $S
\subseteq T$ and $v \in V$ then: 
\[
f(T \cup \{v\}) - f(T) \leq f(S \cup \{v\}) - f(S),
\]
so the effect of each individual is decreasing when the set
increases. 
\begin{assumption}
Throughout, we assume that $f_v(\emptyset) = 0$ and that $f_v$ is monotone
and submodular for all $v \in V$.
\end{assumption}
 
We will consider the following process, which is often referred to as a
``diffusion'' process in the sociology literature. 
\begin{definition}[Diffusion] \label{def:process}
For a given $\fcal$, consider the following process, denoted $\bfss = 
(S_t)_{t=0}^{n-1}$, started
at $S\subseteq V$:
\begin{enumerate}
\item Associate to each node $v$ an independent random
variable  $\theta_v$ uniform in $[0,1]$ ;
\item Set $S_0 = S$;
\item At time $t \geq 1$, initialize $S_t = S_{t-1}$ and add to $S_{t}$ the set of nodes
in $V\sm S_{t-1}$ such that $f_v(S_{t-1}) \geq \theta_v$.
\end{enumerate}
Clearly the process stops on or before time $n-1$.
We denote by $\qcal_\fcal(S)$ the {distribution of} 
$\bfss$ when started at $S$ and
write $\bfss \sim \qcal_\fcal(S)$, where we will drop
the subscript when $\fcal$ is clear from the context.
\end{definition}

\begin{definition}[Influence]
For a weight function $w:2^V \to \real_+$, 
we define the {\it influence} $\sigma_w(S)$ of $S\subseteq V$
as
\begin{equation*}
\sigma_w(S) = \expec_S[w(S_{n-1})],
\end{equation*} 
where $\expec_S$ is the expectation under $\qcal_\fcal(S)$.
\end{definition}

\subsection{Previous Results}

Conjecture~\ref{conj} was previously verified in several special
cases and related models.

{\bf Linear Threshold Model~\cite{KeKlTa:03}.} This is the general
threshold model with $f_v$ of the form
\begin{equation*}
f_v(S) = \sum_{w\in S} b_{v,w},
\end{equation*}
for nonnegative constants $b_{v,w}$. The proof uses a representation
in terms of a related percolation model. See~\cite{KeKlTa:03} for details. 

{\bf ``Normalized'' Submodular Threshold Model~\cite{KeKlTa:05}.} This
is the general threshold model with $f_v$ satisfying the so-called
``normalized'' submodularity property:
\begin{equation}\label{eq:normalized}
\frac{f_v(S \cup \{i\}) - f_v (S)}{1 - f_v(S)}
\geq
\frac{f_v(T \cup \{i\}) - f_v (T)}{1 - f_v(T)},
\end{equation}
for all $S\subseteq T$. Note that this is stronger
than submodularity. The proof takes advantage of an equivalence
with the {\em decreasing cascade model} (see below).

{\bf Independent Cascade Model~\cite{KeKlTa:03}.}
This is a related model where each edge $(v,w)$ has an associated
probability $p_{v,w}$ of being {\em live}, independently of all other edges.
Infected nodes are those connected to the initial set through a 
{\em live path}. The proof of Conjecture~\ref{conj} in this
case also uses a percolation argument.

{\bf Decreasing Cascade Model~\cite{KeKlTa:05}.}
A natural generalization of the previous model consists in defining
for each $v$, each neighbor $w$ of $v$ and each 
subset of neighbors $S$ of $v$ a success probability
$p_v(w, S)$ which is the probability that node $w$ will
succeed in activating $v$ given that nodes in $S$ are active
and have failed to activate $v$. Each node $w$ gets only 
one chance to activate each of its neighbors. 
KKT impose a natural {\em order-independence} condition
on the success probabilities, that is, the overall success
probability of activating $v$ does not depend on the order in which the
active neighbors of $v$ try to activate it. This model---called
the {\em general cascade model} in~\cite{KeKlTa:03}---turns out 
to be equivalent to the {\em general threshold model} under the maps
\begin{equation*}
p_v(w, S) = \frac{f_v(S \cup \{w\}) - f_v (S)}{1 - f_v(S)},
\end{equation*}
and
\begin{equation*}
f_v(S) = 1 - \prod_{i=1}^r (1 - p_v(w_i, S_{i-1})),
\end{equation*}
where $S = \{w_1,\ldots,w_r\}$ and $S_i = \{w_1,\ldots,w_i\}$.
When
\begin{equation}\label{eq:decreasing}
p_v(w,S) \geq p_v(w,T)
\end{equation} 
for all $S\subseteq T$ and all $v,w$, 
the model is called
the {\em decreasing cascade model}. 

It is easy to check
that the decreasing cascade model is equivalent to (\ref{eq:normalized}) 
under the mapping above.
The proof of the conjecture for the decreasing cascade model
works by coupling the processes started at $S$ and $T$ with
$S\subseteq T$ and then adding $w$ in a second phase where condition
(\ref{eq:decreasing}) is used.

% El: Changed wording here. 

In~\cite{KeKlTa:03}, it is also shown that these results carry over to the
{\em non-progressive} case where $\theta_v$ is resampled independently 
at each time step and to {\em general marketing strategies} where
one can use several marketing actions simultaneously. We refer the reader to~\cite{KeKlTa:03}
for details.

{\bf Contact Process.} 
A result similar to our Theorem~\ref{thm:main} below 
holds for a related {\em contact process}~\cite{Harris:74, Liggett:85,LiScSc:u},
where the infection rates are monotone, concave functions of the number of
infected neighbours (a special case of monotone, submodular functions) 
and where vertices also heal at a constant rate.
In particular, Harris~\cite{Harris:74} uses 
a coupling argument, although it appears that
Harris' proof does not extend easily to our setting
as it uses the special form of the infection rates.

%\vfill\eject

{\bf Subsequent Work.} After the publication of this work
in its extended abstract form~\cite{MosselRoch:07b}, further
work was devoted to viral marketing, particularly in the
competitive framework where firms compete
for customers~\cite{DuGadMe:06,CaNaWivZu:07,BhKeSa:07,Even-DarShapira:07}.

\subsection{Main Result}

\begin{theorem}[Main Result]\label{thm:main}
Consider the process defined in Definition~\ref{def:process} 
where $\fcal$ and $w$ are monotone and submodular.
Then, $\sigma_w$ is monotone and submodular.
In particular, this is true when $w$ is the cardinality function.
\end{theorem}

\begin{corollary}\label{cor}
Consider the process defined in Definition~\ref{def:process} 
where $\fcal$ and $w$ are monotone and submodular. 
Furthermore, assume that $w$ takes values between
$1$ and $\poly(n)$.
Then, there exists a (greedy) 
$(1-1/e-1/\poly(n))$-approximation algorithm for maximizing 
$\sigma_w(S)$ among all sets $S$ of size $k$~\cite{KeKlTa:05}
in time $\poly(n)$. 
In particular, this is true when $w$ is the cardinality function.
\end{corollary}
The corollary follows from Theorem~\ref{thm:main} and Theorem 2
of~\cite{KeKlTa:05}. KKT's Greedy Approximation Algorithm is a simple variant
of the standard greedy algorithm where sampling is used to
estimate $\sigma_w$. 

{

{\bf Our proof.} Similarly to~\cite{KeKlTa:05}, a natural idea is to run the process in 
stages. Here we use three phases: we first grow $A\cap B$, then 
$A\setminus B$, and finally $B\setminus A$. See Figure~\ref{fig:coupling} 
for an illustration. The key difference is in the execution of the last 
phase. To do away with the ``normalized'' submodularity condition 
of~\cite{KeKlTa:05}, we do the following.
\begin{itemize}

\item We use a careful combination of cascade and threshold models, which we call 
the need-to-know representation.

\item More importantly, we introduce a novel ``antisense'' coupling technique based on the 
intuition that coupling the processes started at arbitrary sets $A$ and 
$B$ by using $\theta_v$ and $1-\theta_v$ respectively, in a way, 
``maximizes their union'' (note that $1-\theta_v$ is also uniform in 
$[0,1]$). This has to be implemented carefully to also control the 
intersection. See Section~\ref{sec:proof} for details.
See, e.g.,~\cite{Lindvall:92} for a general reference on the
coupling method.
\end{itemize}
}

The rest of the paper is organized as follows. We begin with a few preliminary remarks
in Section~\ref{sec:piecemeal}. The need-to-know representation
and the antisense process are introduced 
in Section~\ref{sec:antiphase}. The full coupling and the proof of its 
correctness appear in Section~\ref{sec:coupling}. 

The results
of this paper were announced in the form of an extended abstract in~\cite{MosselRoch:07b}.

\section{Proof}\label{sec:proof}

Throughout we fix $\fcal$ and $w$ monotone, submodular.
We also fix two arbitrary sets $A, B \subseteq V$ and let $C = A\cap B$
and $D = A\cup B$. 
The idea of the proof is to couple the four processes
\begin{equation*}
\begin{array}{llll}
\bfaa = (A_t)_{t=0}^{n-1} \sim \qcal(A),\\
\bfbb = (B_t)_{t=0}^{n-1} \sim \qcal(B),\\
\bfcc = (C_t)_{t=0}^{n-1} \sim \qcal(C),\\
\bfdd = (D_t)_{t=0}^{n-1} \sim \qcal(D),
\end{array}
\end{equation*}
in such a way that 
\begin{equation}\label{eq:condinter}
C_{n-1} \subseteq A_{n-1}\cap B_{n-1},
\end{equation}
and
\begin{equation}\label{eq:condunion}
D_{n-1} \subseteq A_{n-1} \cup B_{n-1}.
\end{equation}
Indeed, we then have the following lemma.
\begin{lemma}\label{lem:finalstep}
Suppose there exists a coupling of $\bfaa,\bfbb,\bfcc$ and $\bfdd$ 
satisfying (\ref{eq:condinter}), (\ref{eq:condunion}). Then 
\begin{equation} \label{eq:sigma}
\sigma_w(A) + \sigma_w(B) \geq \sigma_w(A\cap B) + \sigma_w(A\cup B).
\end{equation}
\end{lemma}
\begin{proof}
We have by monotonicity and submodularity of $w$
\begin{eqnarray}
&&w(A_{n-1}) + w(B_{n-1})\nonumber\\
&&\qquad\qquad\geq w(A_{n-1}\cap B_{n-1}) + w(A_{n-1}\cup B_{n-1})\nonumber\\
&&\qquad\qquad\geq w(C_{n-1}) + w(D_{n-1}),
\end{eqnarray}
and therefore, taking expectations we get (\ref{eq:sigma}).
\end{proof}

Our coupling is based on the following ideas:
\begin{itemize}
\item[-] {\bf Antisense coupling.} The obvious coupling is to use the same $\theta_v$'s for all processes. It is easy 
to see that such a coupling does not satisfy (\ref{eq:condunion}). It does however
satisfy (\ref{eq:condinter}). 
Intuitively, using the same $\theta_v$ for $\bfaa$ and $\bfbb$ ``maximizes
their intersection'' while using $\theta_v$ for $\bfaa$ and $(1 - \theta_v)$ for
$\bfbb$ ``maximizes their union.'' 
We call this last coupling the {\it antisense coupling}.
To dominate both the intersection and the union simultaneously,
we combine these two couplings. 

\item[-] {\bf Piecemeal growth.} The growth of the four processes can be divided in several stages where we add
the initial sets progressively. Roughly, 
the coupling below starts by growing $A\cap B$, then $A\sm B$ and finally
$B\sm A$. Following our previous comment, the last phase uses 
the antisense coupling to allow the process $\bfbb$ to dominate $\bfdd$ in that
phase.

\item[-] {\bf Need-to-know representation.} Finally, to help carry out the  previous remarks, 
we note that it is not necessary to pick the
$\theta_v$'s at the beginning of the process. Instead, at each step, we uncover as little information
as possible about $\theta_v$. {This is related to the cascade model of~\cite{KeKlTa:05}
although here we use an explicit combination of cascade and threshold models.} 
\end{itemize}
We explain these ideas next. The proof of Theorem~\ref{thm:main} follows in Section~\ref{sec:coupling}.

\subsection{Piecemeal growth}\label{sec:piecemeal}

We first describe an equivalent representation of the process where the initial
set is added in stages. 
We denote by $\qcal(S\,|\,\theta)$ the process $\qcal(S)$ conditioned on
$\theta = (\theta_v)_{v\in V}$.
For a partition $S^{(1)},\ldots,S^{(K)}$ of $S$ (we allow some of the $S^{(k)}$'s to be empty),
consider the process
\begin{equation*}
\bftt =(T_t)_{t=0}^{Kn - 1} \sim \qcal(S^{(1)}, \ldots, S^{(K)}),
\end{equation*}
where
\begin{enumerate}
\item For each $v\in V$ pick $\theta_v$ uniformly and independently in $[0,1]$ and set $T_{-1} = \emptyset$;
\item For $1 \leq k\leq K$, we set 
\begin{equation*}
(T_t)_{t=(k-1)n}^{kn - 1} \sim \qcal(T_{(k-1)n - 1}\cup S^{(k)}\,|\,\theta);
\end{equation*} 
in other words, we add the $S^{(k)}$'s one at a time and use the same $\theta_v$'s for all stages.
\end{enumerate}

The outcomes of the processes $\qcal(S)$ and
$\qcal(S^{(1)},\ldots,S^{(K)})$ have the same distribution.
This result actually follows from a more general discussion
in~\cite{KeKlTa:05}, but we give a proof here for
completeness.
\begin{lemma}[Piecemeal Growth]\label{lem:piecemeal}
Let $S^{(1)},\ldots,S^{(K)}$ be a partition of $S \subseteq V$.
Let
\begin{equation*}
\bfss=(S_t)_{t=0}^{n - 1} \sim \qcal(S),
\end{equation*}
and
\begin{equation*}
\bftt=(T_t)_{t=0}^{Kn - 1} \sim \qcal(S^{(1)}, \ldots, S^{(K)}).
\end{equation*}
Then $S_{n-1}$ and $T_{Kn - 1}$ {have the same distribution}.
\end{lemma}
\begin{proof}
Pick $\theta_v$ uniformly and independently in $[0,1]$ for each $v\in V$ and
let 
\begin{equation*}
\bfss=(S_t)_{t=0}^{n - 1} \sim \qcal(S\,|\,\theta),
\end{equation*}
and
\begin{equation*}
\bftt=(T_t)_{t=0}^{Kn - 1} \sim \qcal(S^{(1)}, \ldots, S^{(K)}\,|\,\theta).
\end{equation*}
Moreover, let
\begin{equation*}
\bftt'=(T'_t)_{t=0}^{Kn - 1} \sim \qcal(S, \emptyset,\ldots,  \emptyset\,|\,\theta),
\end{equation*}
and
\begin{equation*}
\bftt''=(T''_t)_{t=0}^{Kn - 1} \sim \qcal( \emptyset,\ldots,\emptyset, S\,|\,\theta).
\end{equation*}
By monotonicity and induction on the $K$ stages,
\[
T''_{Kn-1}\subseteq T_{Kn-1}\subseteq T'_{Kn-1}
\] 
But clearly 
\[
T'_{Kn-1} = T''_{Kn-1} = S_{n-1}
\] 
so that $S_{n-1} = T_{Kn-1}$.
\end{proof}

\subsection{Antisense phase and need-to-know representation}\label{sec:antiphase}

To implement the antisense coupling, we define the following
variant of the process. 
%See Figure~\ref{fig:antisense} for an illustration.
\begin{definition}\label{def:antisense}
Let $S^{(1)},\ldots,S^{(K)}$ be a partition of $S$ and let $T \subseteq V\sm S$.
We define the process 
\[
\bftt = (T_t)_{t=0}^{(K+1)n-1} \sim \qcal_{-}(S^{(1)},\ldots,S^{(K)}; T)
,
\]
where 
\begin{enumerate}
\item For each $v\in V$ pick $\theta_v$ uniformly in $[0,1]$;

\item Let $\bftt=(T_t)_{t=0}^{Kn - 1} \sim \qcal(S^{(1)}, \ldots, S^{(K)}\,|\,\theta)$;

\item Set $T_{Kn} = T_{Kn - 1} \cup T$;

\item At time $Kn + 1\leq t \leq (K+1)n - 1$, initialize $T_t = T_{t-1}$, and
add to $T_{t}$ the set of nodes
in $V\sm T_{t-1}$ such that 
\[
f_v(T_{t-1}) - f_v(T_{Kn-1}) \geq 1 - \theta_v.
\]

\end{enumerate}
\end{definition}
%\begin{figure}
%\centering
%\includegraphics[width=0.75\textwidth]{antisense.jpg}
%\caption{Schematic representation of the antisense phase. Top: original process seen from a node $v$;
%the activation function gradually fills up the $[0,1]$ interval. Bottom: antisense process seen from a node
%$v$; in the last phase, the interval is filled up ``from the end.''}\label{fig:antisense}
%\end{figure}

\begin{lemma}[Antisense Phase]\label{lem:antiphase}
Assume $S^{(1)},\ldots,S^{(K)}$ is a partition of $S$ and
$T \subseteq V\sm S$.
Let 
\begin{equation*}
\bfss = (S_t)_{t=0}^{(K+1)n -1} 
\sim \qcal(S^{(1)},\ldots,S^{(K)}, T),
\end{equation*}
and
\begin{equation*}
\bftt = (T_t)_{t=0}^{(K+1)n-1}
\sim \qcal_{-}(S^{(1)},\ldots,S^{(K)}; T).
\end{equation*}
Then, $S_{(K+1)n - 1}$ and $T_{(K+1)n - 1}$ {have the same distribution}.
\end{lemma}
\begin{proof}
As was discussed at the beginning of Section~\ref{sec:proof}, 
rather than picking the $\theta_v$'s at the beginning of the process,
it is useful to think of them as being progressively uncovered on
a need-to-know basis. Consider only the {\em first} stage of the 
process $\bfss$ for the time being. Let $S_{-1} = \emptyset$.
Suppose that, at time $t \geq 1$, $v \notin S_{t-1}$. Then
we have that $\theta_v \in [f_v(S_{t-2}), 1]$ and all we need to
know to decide if $v$ is added to $S_t$ is whether or not
$\theta_v \in [f_v(S_{t-2}), f_v(S_{t-1})]$. In other words,
was the increase in $f_v$ between time $t-2$ and $t-1$ enough
to hit $\theta_v$? Note that, given the event $\{ f_v(S_{t-2}) \leq \theta_v\}$,
$\theta_v$ is uniformly distributed in $[f_v(S_{t-2}),1]$ and
we have that $\theta_v$ is in $[f_v(S_{t-2}), f_v(S_{t-1})]$ with probability
\begin{equation*}
\frac{f_v(S_{t-1}) - f_v(S_{t-2})}{1 - f_v(S_{t-2})}.
\end{equation*}
Therefore, we can describe the process $(S_t)_{t=0}^{n - 1}$ equivalently
as follows.
We first set $S_{-1} = \emptyset$, $S_0 = S$.
Then, at step $1 \leq t\leq n-1$, 
we initialize $S_t = S_{t-1}$ and 
for each $v \in V \sm S_{t-1}$:
\begin{itemize}
\item[-] With probability 
\begin{equation} \label{eq:fancy_proc1}
\frac{f_v(S_{t-1}) - f_v(S_{t-2})}{1-f_v(S_{t-2})},
\end{equation}
we add $v$ to $S_t$ and pick $\theta_v$ uniformly
in
\begin{equation*}
[f_v(S_{t-2}), f_v(S_{t-1})];
\end{equation*} 
\item[-] Otherwise, we do nothing.
\end{itemize}
By the discussion above, this new version of the process {has the 
same distribution as $\qcal(S^{(1)})$.} 
We proceed similarly for the following $K-1$ stages to get
$(S_t)_{t=0}^{Kn-1}$ which is then distributed according to $\qcal(S^{(1)},\ldots,S^{(K)})$.

Up to time $Kn-1$, the processes $(S_t)$ and $(T_t)$
have identical transition probabilities. Hence, we can take 
\[
(T_t)_{t=0}^{Kn-1} = (S_t)_{t=0}^{Kn-1}.
\]
Then note that, at time $t = Kn$, for each $v\notin S_{Kn-1}=T_{Kn-1}$, 
we have that $\theta_v$ is uniformly distributed in 
\[
[f_v(S_{Kn-1}),1] = [f_v(T_{Kn-1}),1]. 
\]
For each such $v$, we now pick $\theta_v$
uniformly in $[f_v(S_{Kn-1}),1]$. Moreover,
we define for all $v\in V$ 
\begin{equation*}
\theta'_v =
\left\{
\begin{array}{ll}
\theta_v,& v\in S_{Kn-1},\\
f_v(S_{Kn - 1}) + 1-\theta_v,& v\notin S_{Kn-1}. 
\end{array}
\right.
\end{equation*}
Finally, let
\begin{equation*}
(S_t)_{t=Kn}^{(K+1)n - 1} \sim \qcal(S_{Kn - 1}\cup T\,|\,\theta),
\end{equation*}
and
\begin{equation*}
(T_t)_{t=Kn}^{(K+1)n - 1} \sim \qcal(T_{Kn - 1}\cup T\,|\,\theta').
\end{equation*}
That is, we run the last stage of $\bfss$ and $\bftt$ as before,
with $\theta$ and $\theta'$ respectively.
It is clear that $\bftt \sim \qcal_-(S^{(1)},\ldots,S^{(K)};T)$
by construction.
Moreover, it follows that $S_{(K+1)n-1}$ and $T_{(K+1)n-1}$
{have the same distribution}
from the fact that for a uniform variable $\theta_v$
in $[f_v(S_{Kn-1}), 1]$, the random variables $\theta_v$ and $f_v(S_{Kn - 1}) + 1-\theta_v$ 
{have the same distribution.}
\end{proof}

\subsection{Coupling}\label{sec:coupling}

We are now ready to prove Theorem~\ref{thm:main}.
We will need the following easy consequence of monotone submodularity.
\begin{lemma}\label{lem:gensubmod}
Let $f : 2^V \to \real_+$ be monotone and submodular.
Then if $S \subseteq S' \subseteq V$ and $T \subseteq T' \subseteq V$, we have
\begin{equation*}
f(S \cup T') - f(S) \geq f(S' \cup T) - f(S'). 
\end{equation*}
\end{lemma}
\begin{proof}
Note that by monotonicity and submodularity
\begin{eqnarray*}
f(S\cup T') - f(S)
&\geq& f(S\cup T) - f(S)
= f(S\cup (T\sm S)) - f(S)\\
&\geq& f(S \cup (S'\sm (T\cup S)) \cup (T\sm S))
- f(S \cup (S'\sm (T\cup S)))\\
&\geq& f(S'\cup T) - f(S').
\end{eqnarray*}
\end{proof}
We now give a proof of Theorem~\ref{thm:main}.

\begin{proof}
We proceed with our coupling of $\bfaa$, $\bfbb$, $\bfcc$, and $\bfdd$.
In fact, by Lemmas~\ref{lem:finalstep}, \ref{lem:piecemeal}, and~\ref{lem:antiphase}, 
it suffices instead to couple 
\begin{equation*}
\begin{array}{lll}
&&\bfaa = (A_t)_{t=0}^{3n-1}  
\sim \qcal(A\cap B, A \sm B, \emptyset),\\
&&\bfbb = (B_t)_{t=0}^{3n-1} 
\sim \qcal_-(A \cap B, \emptyset; B \sm A),\\
&&\bfcc = (C_t)_{t=0}^{3n-1} 
\sim \qcal(A\cap B,\emptyset, \emptyset),\\
&&\bfdd = (D_t)_{t=0}^{3n-1} 
\sim 
\qcal_-(A\cap B, A\sm B; B\sm A),
\end{array}
\end{equation*}
in such a way that for all $0 \leq t \leq 3n-1$
\begin{eqnarray}\label{eq:submod}
&&C_t \subseteq A_t \cap B_t,
\qquad D_t \subseteq A_t \cup B_t.
\end{eqnarray}
Our coupling is as follows. We pick $\theta_v$ uniformly in $[0,1]$
for all $v\in V$ and use the same $\theta$ for all
four processes above. See Figure~\ref{fig:coupling} for a graphical
representation of the proof.\begin{figure*}
\centering
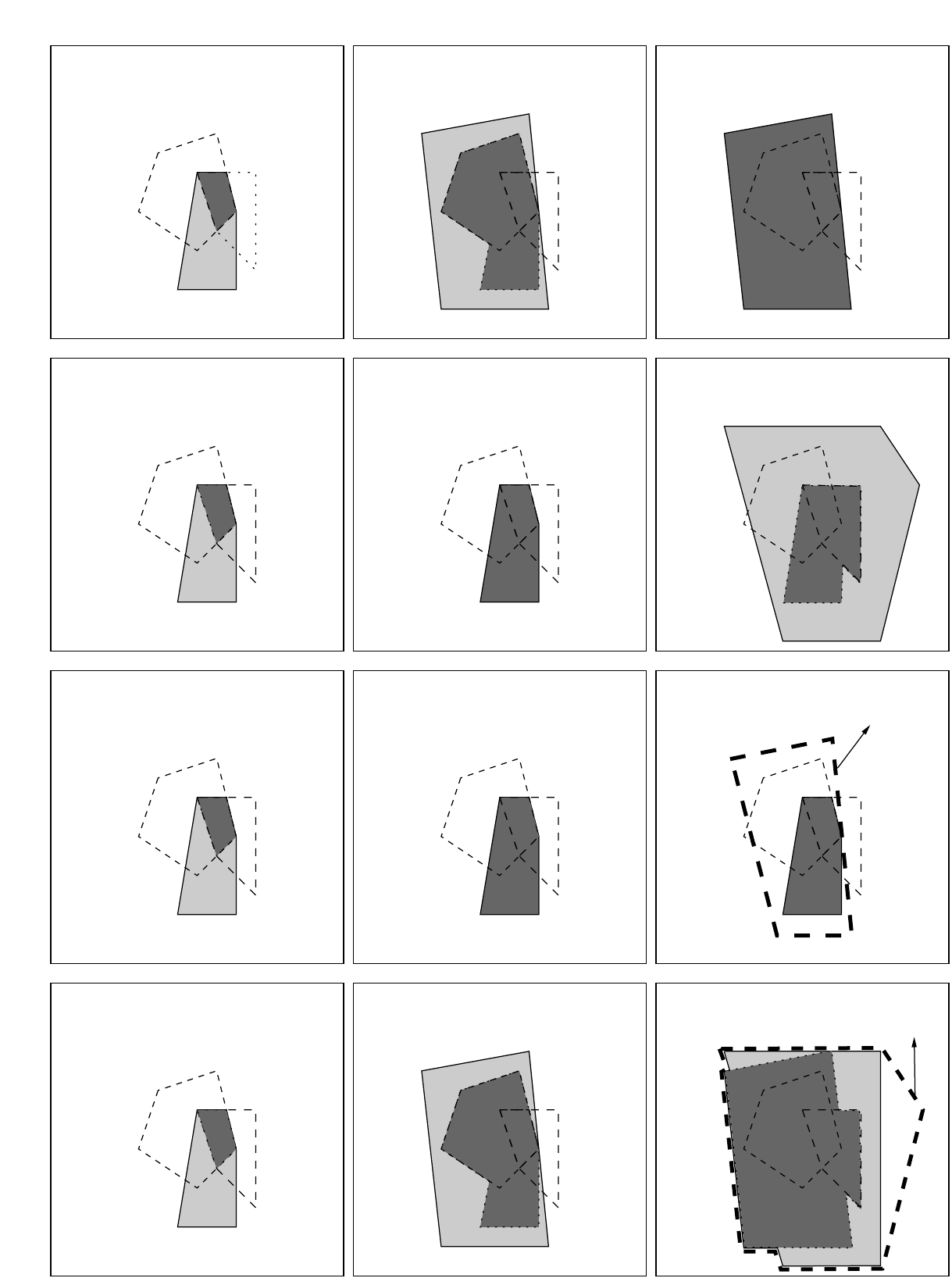\caption{The three phases of the coupling.
In each phase, the dark shaded region is the initial set, while
the light shaded region is the final set. The sets $A$ and $B$
are indicated by dashed lines. {The thick dashed lines
show that the desired properties are satisfied.}}
\label{fig:coupling}
\end{figure*}

By construction, for all $0 \leq t \leq 2n-1$ we have 
\[
B_t = C_t \subseteq A_t
\]
so that 
\[
C_t = A_t \cap B_t. 
\]
Similarly for all $0 \leq t
\leq 2n-1$ we have $D_t = A_t$ so that 
\[
D_t \subseteq A_t \cup B_t.
\]  
Thus (\ref{eq:submod}) is satisfied for $0 \leq t \leq 2n-1$. 
To see (\ref{eq:submod}) holds also for $2n \leq t\leq 3n-1$, 
note that by Lemma~\ref{lem:gensubmod}
for all $v \notin D_{2n}$
\begin{equation*}
f_v(B_{2n}) - f_v(B_{2n-1}) \geq f_v(D_{2n}) - f_v(D_{2n-1}),
\end{equation*}
since $B_{2n-1} \subseteq D_{2n-1}$, 
\[
B_{2n} = B_{2n-1} \cup (B \sm A),
\] and
\[
D_{2n} = D_{2n-1} \cup (B \sm A).
\]
(See Figure~\ref{fig:submodularity} for an illustration of this step.)
We proceed by induction.
By monotonicity and Lemma~\ref{lem:gensubmod}, we then have
for all $2n \leq t\leq 3n-1$
\begin{eqnarray*}
&&(\Omega_{1,t})\qquad D_t\sm D_{2n-1} \subseteq B_t\sm B_{2n-1},\\
%\end{equation*}
%and
%\begin{equation*}
&&(\Omega_{2,t})\qquad f_v(B_{t}) - f_v(B_{2n-1}) \geq f_v(D_{t}) - f_v(D_{2n-1}),
\end{eqnarray*}
for all $v \notin D_{2n}$. 
Indeed, assume $(\Omega_{1,t'})$ and $(\Omega_{2,t'})$ for all
$2n \leq t' \leq t$. We have already proved the base case $t=2n$.
Then $(\Omega_{2,t})$ implies $(\Omega_{1,t+1})$ by definition
of the process in the antisense phase (Definition~\ref{def:antisense}). 
In turn, $(\Omega_{1,t+1})$ implies
$(\Omega_{2,t+1})$ by monotonicity and Lemma~\ref{lem:gensubmod}.

This proves the claim since we then have 
for all $2n \leq t\leq 3n-1$, 
$A_{t} = D_{2n-1}$ and 
\[
D_t\sm D_{2n-1} \subseteq B_t
\]
which implies
\[
D_t \subseteq A_t \cup B_t. 
\]
The condition 
\[
C_t \subseteq A_t\cap B_t
\]
is clear from the construction.
\begin{figure}
\centering
\includegraphics[width=0.75\textwidth]{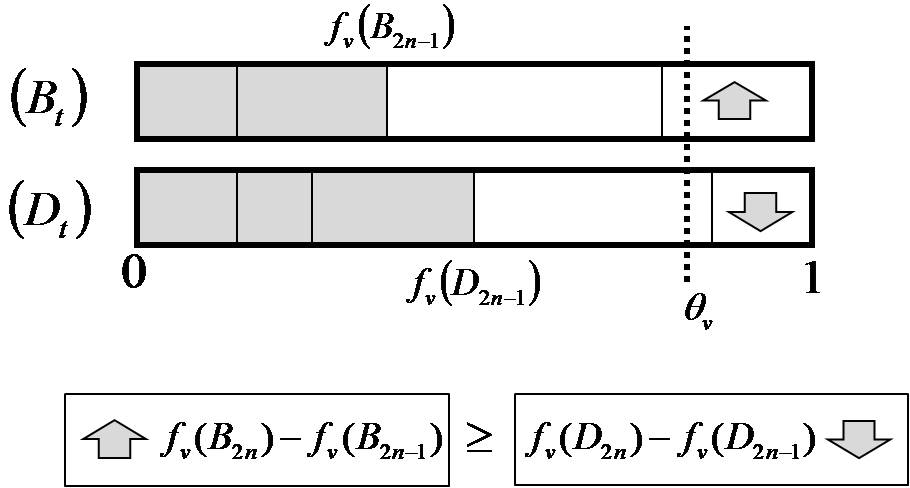}
\caption{Last phase of the coupling seen from a node $v$. 
In the antisense phase, the interval is filled up ``from the end.''
By submodularity, the process $(B_t)$ reaches the
threshold first.}\label{fig:submodularity}
\end{figure}
\end{proof}

\section{Concluding Remarks}\label{sec:conclusion}

{\bf Necessity.}
It is easy to see that the submodularity assumption in 
Theorem~\ref{thm:main} is  necessary in the following sense: 
Any function $f$ which
is not submodular admits a network with activation function $f$ where the influence is not submodular.
Indeed, let $f : 2^V \to \real_+$, $A, B \subseteq V$ such that
\begin{equation*}
f(A) + f(B) < f(A\cap B) + f(A \cup B).
\end{equation*}
Let $V^* = V \cup \{v^*\}$ with $f_{v^*} \equiv f$ and $f_v \equiv 1$ for all
$v \in V$. It is then immediate to check that:
\begin{eqnarray*}
\sigma(A) + \sigma(B) 
&=& |A| + |B| + f(A) + f(B)\\ 
&=& |A \cap B| + |A \cup B|\\ 
&& \qquad + f(A) + f(B)\\ 
&<& |A \cap B| + |A \cup B|\\ 
&& \qquad + f(A \cap B) + f(A \cup B)\\ 
&=& \sigma(A \cap B) + \sigma(A \cup B).
\end{eqnarray*}

\noindent{\bf Other threshold distributions.}
As noted in the introduction,
our results hold for general threshold distributions.
Assume that $\theta_v$ has cumulative distribution
$F_v$. Notice that the dynamics can be re-written as
\begin{equation*}
f_v(S_{t-1}) \geq \theta_v \equiv F_v^{-1}(U_v) \qquad 
\Leftrightarrow \qquad F_v(f_v(S_{t-1})) \geq U_v,
\end{equation*}
where the $U_v$'s are independent uniform in $[0,1]$.
Hence, if $F_v\circ f_v$ is increasing submodular for all $v$
then the influence function is submodular.
Furthermore, the example above demontrates
that this condition is necessary.

\section*{Acknowledgments}

We thank Jon Kleinberg and \'Eva Tardos for insightful discussions
and encouragements. We also thank Tom Liggett for bringing to our attention 
the work of Harris.

\bibliographystyle{alpha}
\bibliography{thesis}

\end{document}